\documentclass{amsart}
\usepackage{amsmath,amssymb}
\newtheorem{theorem}{Theorem}[section]
\newtheorem{lemma}[theorem]{Lemma}
\newtheorem{proposition}[theorem]{Proposition}

\theoremstyle{definition}

\theoremstyle{remark}
\newtheorem{remark}[theorem]{Remark}

\numberwithin{equation}{section}

\begin{document}

\title{On the zeros of the Riemann zeta function}

\author{Lazhar Fekih-Ahmed}
\address{\'{E}cole Nationale d'Ing\'{e}nieurs de Tunis, BP 37, Le Belv\'{e}d\`{e}re 1002 , Tunis, Tunisia}
\curraddr{\'{E}cole Nationale d'Ing\'{e}nieurs de Tunis, BP 37, Le
Belv\'{e}d\`{e}re 1002 , Tunis, Tunisia}
\email{lazhar.fekihahmed@enit.rnu.tn}


\date{January 21, 2011}



\begin{abstract}
This paper is divided into two independent parts. The first part
presents new integral and series representations of the Riemaan
zeta function. An equivalent formulation of the Riemann hypothesis
is given and few results on this formulation are briefly outlined.
The second part exposes a totally different approach. Using the
new series representation of the zeta function of the first part,
exact information on its zeros is provided.
\end{abstract}

\maketitle

\specialsection*{PART I}
\section{Introduction}
\label{intro}

It is well known that the Riemaan zeta function  defined  by the
Dirichlet series
\begin{equation}\label{sec1-eq1} \zeta(s)=
\frac{1}{1^{s}}+\frac{1}{2^{s}}+\cdots+\frac{1}{k^{s}}+\cdots=\sum_{n=1}^{\infty}n^{-s}
\end{equation}
 converges for $\Re(s)>1$, and can be analytically
continued to the whole complex plane with one singularity, a
simple pole with residue $1$ at $s=1$. It is also well known that
$\zeta(s)$ satisfies the functional equation:

\begin{equation}\label{sec1-eq2}
\chi(s)\zeta(s)=\zeta(1-s)\chi(1-s)\quad {\rm with}\quad\chi(s)=
\pi^{-\frac{s}{2}}\Gamma(\frac{s}{2}),
\end{equation}
and that the  zeros of $\zeta(s)$ come into two types. The trivial
zeros which occur at all negative even integers $s=-2,-4,\cdots$,
and the nontrivial zeros which occur at certain values of $s\in
\mathbb{C}$, $0<\Re(s)<1$.

The Riemann hypothesis states that the nontrivial zeros of
$\zeta(s)$ all have real part $\Re(s)=\frac{1}{2}$. From the
functional equation (\ref{sec1-eq2}), the Riemann hypothesis is
equivalent to  $\zeta(s)$ not having any zeros in the strip
$0<\Re(s)<\frac{1}{2}$.

\section{An Analytic Continuation of  $\zeta(s)$}\label{sec2}

Let

\begin{equation}\label{sec2-eq1}
S_{n}(s)=1-{n-1\choose 1}2^{-s} + {n-1\choose 2}3^{-s} -\cdots
+(-1)^{n-1}(n)^{-s}, n\ge 2
\end{equation}

with $S_{1}(s)=1$.

Using the well-known identity, valid for
$\displaystyle{\Re(s)>0}$:

\begin{equation}\label{sec2-eq2}
n^{-s}=\frac{1}{\Gamma(s)}\int_{0}^{\infty}e^{-nt}t^{s-1}\,dt,
\end{equation}

we can rewrite $S_{n}(s)$ in (\ref{sec2-eq1}) as:

\begin{eqnarray}
S_{n}(s) &= &\sum_{k=0}^{n-1}(-1)^{k}{n-1\choose k}(k+1)^{-s} \nonumber\\
   &=& \frac{1}{\Gamma(s)}\int_{0}^{\infty}\sum_{k=0}^{n-1}(-1)^{k}{n-1\choose k}
   e^{-(k+1)t}t^{s-1}\,dt \label{sec2-eq3}\\
  &=&
  \frac{1}{\Gamma(s)}\int_{0}^{\infty}(1-e^{-t})^{n-1}e^{-t}t^{s-1}\,dt,\nonumber
\end{eqnarray}

since  $\displaystyle{\sum_{k=0}^{n-1}(-1)^{k}{n-1\choose k}
e^{-(k+1)t}=e^{-t}(1-e^{-t})^{n-1}}$.

We have
\begin{eqnarray}\label{sec2-eq5}
\sum_{n=1}^{\infty}\frac{S_{n}(s)}{n+1}&=&\frac{1}{\Gamma(s)}\sum_{n=1}^{\infty}\int_{0}^{\infty}\frac{(1-e^{-t})^{n-1}}{n+1}e^{-t}t^{s-1}\,dt\\
&=&\frac{1}{\Gamma(s)}\int_{0}^{\infty}\sum_{n=1}^{\infty}\frac{(1-e^{-t})^{n-1}}{n+1}e^{-t}t^{s-1}\,dt\label{sec2-eq6}\\
&=& \frac{1}{\Gamma(s)}\int_{0}^{\infty}\bigg
(\frac{t}{(1-e^{-t})^2}-\frac{1}{1-e^{-t}}\bigg )
e^{-t}t^{s-1}\,dt.\label{sec2-eq7b}
\end{eqnarray}

Before we proceed further, some remarks are  in order:

\begin{remark} The interchange of the  summation and integration
in equation (\ref{sec2-eq6}) is valid because the series
$\sum_{n=1}^{\infty}\int_{0}^{\infty}\frac{(1-e^{-t})^{n-1}}{n+1}e^{-t}t^{s-1}\,dt$
converges absolutely and uniformly for $0<t<\infty$. To see this,
we show uniform convergence for the  dominating series
$\sum_{n=1}^{\infty}\int_{0}^{\infty}\frac{(1-e^{-t})^{n-1}}{n+1}e^{-t}t^{\sigma-1}\,dt$,
$\sigma=\Re(s)$. Indeed, let $K={\rm
max}((1-e^{-t})^{n-1}e^{-t/2})$, $0<t<\infty$. A straightforward
calculation of the derivative shows that
$K=(1-\frac{1}{2n-1})^{n-1}\frac{1}{\sqrt{2n-1}}$ and  is attained
when $e^{-t}=\frac{1}{2n-1}$. Now, for $n\ge 2$ we have
\begin{eqnarray}
\frac{1}{n+1}\int_{0}^{\infty}(1-e^{-t})^{n-1}e^{-t}t^{\sigma-1}\,dt&=&
\frac{1}{n+1}\int_{0}^{\infty}(1-e^{-t})^{n-1}e^{-t/2}
(e^{-t/2}t^{\sigma-1})\,dt\nonumber\\
&\le & \frac{K}{n+1} \int_{0}^{\infty}e^{-t/2}t^{\sigma-1}\,dt\nonumber\\
&=&\frac{1}{n+1} (1-\frac{1}{2n-1})^{n-1}\frac{2^{\sigma}\Gamma(\sigma)}{\sqrt{2n-1}}\label{sec2-eq9}\\
&\le&\frac{2^{\sigma}\Gamma(\sigma)}{(n+1)\sqrt{2n-1}}.\nonumber
\end{eqnarray}

The last inequality implies that the dominating series converges
by the comparison test.
\end{remark}

\begin{remark}
To get equation (\ref{sec2-eq7b}) we used the identity:
\begin{equation}\label{sec2-eq8}
\sum_{n=1}^{\infty}\frac{(1-e^{-t})^{n-1}}{n+1}=\frac{t}{(1-e^{-t})^2}-\frac{1}{1-e^{-t}},
\end{equation}
which can be obtained by putting $X=1-e^{-t}$ into
$\displaystyle{-\frac{\log(1-X)}{X^2}-\frac{1}{X}=\sum_{n=1}^{\infty}\frac{X^{n-1}}{n+1}}$.
\end{remark}

Now, since
\begin{equation}\label{sec2-eq9b}
\frac{d}{dt}\frac{-te^{-t}}{1-e^{-t}}=\frac{te^{-t}}{(1-e^{-t})^2}-\frac{e^{-t}}{1-e^{-t}},
\end{equation}

an integration by parts in (\ref{sec2-eq7b}) yields
\begin{equation}\label{sec2-eq10}
\sum_{n=1}^{\infty}\frac{S_{n}(s)}{n+1} =
\frac{s-1}{\Gamma(s)}\int_{0}^{\infty}\frac{e^{-t}t^{s-1}}{1-e^{-t}}\,dt=(s-1)\zeta(s),
\end{equation}

which is valid when $\Re(s)>1$. And since the integral
(\ref{sec2-eq7b}) is valid for $\Re(s)>0$, then we have proved
\begin{theorem}\label{sec2-thm1}
Let $ \phi(t)=\frac{t}{(1-e^{-t})^2}-\frac{1}{1-e^{-t}}$, then for
$\Re(s)>0$,
\begin{equation}\label{sec2-eq10bb}
(s-1)\zeta(s) =\frac{1}{\Gamma(s)}\int_{0}^{\infty}\phi(t)
e^{-t}t^{s-1}\,dt.
\end{equation}
\end{theorem}

\begin{remark}
Although we will not need it in the rest of the paper, we can also
obtain an analytic continuation of $(s-1)\zeta(s)$ when $\Re(s)\le
0$. We simply rewrite (\ref{sec2-eq10bb}) as a contour integral
\begin{equation}\label{sec2-eq11}
\frac{\Gamma(1-s)}{2\pi
i}\int_{\mathcal{C}}\phi(t)e^{-t}t^{s-1}\,dt,
\end{equation}

where $\mathcal{C}$ is the Hankel contour consisting of the three
parts $C=C_{-}\cup C_{\epsilon}\cup C_{+}$:  a path  which extends
from $(-\infty,-\epsilon)$, around the origin counter clockwise on
a circle of center the origin and of radius $\epsilon$ and back to
$(-\epsilon,-\infty)$, where $\epsilon$ is an arbitrarily small
positive number. The integral (\ref{sec2-eq11}) now defines
$(s-1)\zeta(s)$ for all $s\in\mathbb{C}$.
\end{remark}

\begin{remark}
In particular, when $s=k$ is a positive integer, we have yet
another formula for $\zeta(k)$:

\begin{equation}\label{sec2-eq10b}
(k-1)\zeta(k)=\frac{1}{(k-1)!}\int_{0}^{\infty}\phi(t)e^{-t}t^{k-1}\,dt.
\end{equation}

\end{remark}

\begin{remark} The above integral formula for $(s-1)\zeta(s)$,
although obtained by elementary means, does not seem to be found
in the literature. As for the series formula, it has been obtained
by a different method in \cite{ser}. A series formula that is
different but similar in form and often mentioned in the
literature is that of Hasse \cite{hasse}.
\end{remark}

\section{A Series Expansion of $(s-1)\zeta(s)\Gamma(s)$}\label{sec2bis}

The analytic  function $(s-1)\zeta(s)\Gamma(s)$ can be represented
by a Taylor series around any point $s_0=1+iy$ on the vertical
line $\sigma=1$. In particular, for  $s_0=1$ we obtain the
well-known power series \cite{apostol}:

\begin{equation}\label{sec2bis-eq1}
(s-1)\zeta(s)\Gamma(s)=a_0+a_1(s-1)+a_2(s-1)^2+a_3(s-1)^3+\cdots
\end{equation}

where the coefficients $a_0=1$ and $a_n$ are defined by

\begin{eqnarray}\label{sec2bis-eq2}
a_n&=&\frac{1}{n!}\lim_{s\to 1}\frac{d^n}{d s^n}\bigg \{\int_{0}^{\infty}\phi(t) e^{-t}t^{s-1}\,dt \bigg \}\\
&=&\frac{1}{n!} \int_{0}^{\infty}\phi(t)e^{-t}\lim_{s\to
1}\frac{d^n}{d s^n} \big \{ t^{s-1}\big \}\,dt\\
&=&\frac{1}{n!}\int_{0}^{\infty}\phi(t)e^{-t} (\log{t})^{n}\,dt,
\end{eqnarray}

with $\phi(t)$ being the function defined in
Theorem~\ref{sec2-thm1}.

The  coefficients $a_n$ are very important in the evaluation of
$\zeta^{(k)}(0)$ as given by Apostol in \cite{apostol}. Up to now
the $a_n$ are regarded as unknowns and as difficult to approximate
as $\zeta^{(k)}(0)$ itself as pointed out by Lehmer \cite{lehmer}.
The formula above solves the exact evaluation problem of the
$\zeta^{(k)}(0)$ and many other variant formulae. The following
proposition provides more information on the sequence $\{a_n\}$.

\begin{proposition}\label{sec2bis-prop1}
For $n$ large enough, the coefficients $a_n$ are given by
\begin{equation}\nonumber
a_n = (-1)^n\bigg (\frac{1}{2}-\frac{1}{6}\frac{1}{2^{n+1}}\bigg
)+O(\frac{1}{4^{n}}).
\end{equation}
\end{proposition}
\begin{proof}

The expression of $a_n$ can be split into the sum

\begin{eqnarray}
a_n&=&\frac{1}{n!}\int_{0}^{1}\phi(t)e^{-t}
(\log{t})^{n}\,dt+\frac{1}{n!}\int_{1}^{\infty}\phi(t)e^{-t}
(\log{t})^{n}\,dt\nonumber\\
&=&\frac{(-1)^{n}}{2}+\frac{(-1)^{n}}{n!}\int_{0}^{1}\big
[\phi(t)e^{-t}-\frac{1}{2} \big ] \log \bigg (\frac{1}{t}\bigg
)^{n}\,dt\label{sec2bis-eq3}\\
&&\qquad\qquad +\frac{1}{n!}\int_{1}^{\infty}\phi(t)e^{-t}
(\log{t})^{n}\,dt.\nonumber
\end{eqnarray}

To obtain an estimate the first integral in (\ref{sec2bis-eq3}),
we use equation (\ref{sec2-eq9b}). A differentiation  with respect
to $t$ of the following expansion which defines the Bernoulli
numbers\footnote{$B_0=1$ $B_1=-1/2$, $B_2=1/6$, $B_3=0$,
$B_4=-1/30$, $B_5=0$, $B_6=1/42$, $B_7=0$, $B_8=-1/30$ etc.}

\begin{equation}\label{sec2bis-eq4}
\frac{te^{-t}}{1-e^{-t}}=\sum_{n=0}^{\infty}\frac{B_n}{n!}t^n,\quad
|t|<2\pi
\end{equation}

gives

\begin{equation}\label{sec2bis-eq5}
\phi(t)e^{-t}-\frac{1}{2}=\sum_{n=2}^{\infty}\frac{-B_n}{(n-1)!}t^{n-1}=
-\frac{1}{6}t+\frac{1}{180}t^3-\frac{1}{5040}t^5+\cdots
\end{equation}

Now, since

\begin{equation}\label{sec2bis-eq6}
\frac{1}{n!}\int_{0}^{1} t^m \log \bigg (\frac{1}{t}\bigg
)^{n}\,dt= \frac{1}{(m+1)^{n+1}}
\end{equation}

for all $n, m$ positive, the first integral in (\ref{sec2bis-eq3})
without the factor $(-1)^n$ has an expansion

\begin{eqnarray}\nonumber
\frac{1}{n!}\int_{0}^{1}\big [\phi(t)e^{-t}-\frac{1}{2} \big ]
\log \bigg (\frac{1}{t}\bigg )^{n}\,dt
&=&\sum_{n=2}^{\infty}\frac{-B_n}{(n-1)!(n+1)^{n+1}}\\
&=& -\frac{1}{6}\frac{1}{2^{n+1}}+\frac{1}{180}\frac{1}{4^{n+1}}
-\cdots\label{sec2bis-eq7}
\end{eqnarray}

To estimate the second integral in (\ref{sec2bis-eq3}), we use the
bound

\begin{equation}\label{sec2bis-eq8}
 (\log{t})^n < e^{\epsilon t}
\end{equation}

which is valid for all $t\ge n^{1+\epsilon}$ and for $n$ large
enough and where $\epsilon$ is any positive small number. We split
the integral into two parts

\begin{eqnarray}\label{sec2bis-eq9}
\frac{1}{n!}\int_{1}^{\infty}\phi(t)e^{-t} (\log{t})^{n}\,dt&=&
\frac{1}{n!}\int_{1}^{n^{1+\epsilon}}\phi(t)e^{-t}
(\log{t})^{n}\,dt + \nonumber\\
&&\qquad\qquad\qquad
\frac{1}{n!}\int_{n^{1+\epsilon}}^{\infty}\phi(t)e^{-t} (\log{t})^{n}\,dt\\
&\le&
\frac{C_0}{n!}(1+\epsilon)^n\log(n)^n+\frac{1}{2n!}\frac{e^{-(1-\epsilon)n^{1+\epsilon}}}{1-\epsilon},
\label{sec2bis-eq10}
\end{eqnarray}

where $C_0=\int_{1}^{\infty}\phi(t)e^{-t}\,dt=0.58$.

Clearly, the term
$\frac{1}{2n!}\frac{e^{-(1-\epsilon)n^{1+\epsilon}}}{1-\epsilon}$
is extremely small for $n$ large enough. In particular, for $n\ge
2$, it is less than $\frac{C}{a^n}$ where $C$ is a positive
constant and $a$ is any positive constant greater than 2.

Using Stirling formula $n!\sim \sqrt{2\pi n}\big (\frac{n}{e}\big
)^n$, we can  verify that the term
$\frac{C_0}{n!}(1+\epsilon)^n\log(n)^n$ can also be made less than
$\frac{C}{a^n}$ for $n$ large enough. Taking $a=4$ for example, we
obtain

\begin{equation}\label{sec2bis-eq11}
\frac{1}{n!}\int_{1}^{\infty}\phi(t)e^{-t} (\log{t})^{n}\,dt\le
\frac{C}{4^n}.
\end{equation}

By combining the above estimates, we obtain

\begin{equation}\label{sec2bis-eq12}
a_n = (-1)^n\bigg (\frac{1}{2}-\frac{1}{6}\frac{1}{2^{n+1}}\bigg
)+O(\frac{1}{4^{n}}).
\end{equation}
\end{proof}

For any  $s_0$ of the form $s_0=1+iy$ the corresponding Taylor
series is

\begin{equation}\label{sec2bis-eq13}
(s-1)\zeta(s)\Gamma(s)=b_0+b_1(s-s_0)+b_2(s-s_0)^2+b_3(s-s_0)^3+\cdots
\end{equation}

where the  $b_n$ is expressed as

\begin{equation}\label{sec2bis-eq14}
b_n=\frac{1}{n!}\int_{0}^{\infty}\phi(t)e^{-t}
(\log{t})^{n}t^{iy}\,dt.
\end{equation}

An asymptotic estimate of the coefficients $b_n$ is given by the
following proposition whose proof is the same as the proof of
Proposition~\ref{sec2bis-prop1}.

\begin{proposition}\label{sec2bis-prop2}
The radius of convergence of the series (\ref{sec2bis-eq13}) is
$\sqrt{1+y^2}$ and for $n$ large enough, the coefficients $b_n$
are given by
\begin{equation}\nonumber
b_n =
\frac{1}{2}\frac{(-1)^n}{(1+iy)^{n+1}}-\frac{1}{6}\frac{(-1)^n}{(2+iy)^{n+1}}+O(\frac{1}{(\sqrt{16+y^2})^n}).
\end{equation}
\end{proposition}

\section{The zeros of $\zeta(s)$}\label{sec3bis}

The Taylor series expansion  $(s-1)\zeta(s)\Gamma(s)$ provides us
with a tool to study the zeros of $\zeta(s)$ is a neighborhood of
$s_0=1+iy$. We have
\begin{equation}\label{sec3bis-eq1}
(s-1)\zeta(s)\Gamma(s)=b_0+b_1(s-s_0)+b_2(s-s_0)^2+b_3(s-s_0)^3+\cdots
\end{equation}

with

\begin{equation}\label{sec3bis-eq2}
b_0=iy\zeta(1+iy)\Gamma(1+iy),
\end{equation}

\begin{equation}\label{sec3bis-eq3}
b_n=\frac{1}{n!}\int_{0}^{\infty}\phi(t)e^{-t}
(\log{t})^{n}t^{iy}\,dt.
\end{equation}

It is a well-know fact \cite{hadamard1} that $\zeta(1+iy)\neq 0$
for all $y$. Therefore $b_0\neq 0$ for all $y$ and the inverse of
$(s-1)\zeta(s)\Gamma(s)$ is well-defined and can be expanded into
a power series of the form

\begin{equation}\label{sec3bis-eq4}
\frac{1}{(s-1)\zeta(s)\Gamma(s)}=c_0+c_1(s-s_0)+c_2(s-s_0)^2+c_3(s-s_0)^3+\cdots,
\end{equation}

where the coefficients $c_n$ are given by

\begin{equation}\label{sec3bis-eq5}
c_n=(-1)^n\frac{\Delta_n}{b_0^{n+1}},
\end{equation}
with

\begin{equation}\label{sec3bis-eq6}
\Delta_n=\left|%
\begin{array}{ccccc}
  b_1 & b_0 & \cdots & \cdots & 0 \\
  b_2 & b_1 & b_0 & \cdots & 0\\
  \vdots & \vdots & \vdots & \vdots &\vdots \\
  b_{n-1} & b_{n-2} & b_{n-3} & \cdots & b_{0} \\
  b_{n} & b_{n-1} & b_{n-2} & \cdots & b_{1} \\
\end{array}%
\right|.
\end{equation}

Let $D(s_0,\frac{1}{2})$ be the open disk of center $s_0$ and
radius $\frac{1}{2}$. The zeros of $\zeta(s)$
 are the same as those of  $(s-1)\zeta(s)\Gamma(s)$ in the right
half plane; therefore,  $\zeta(s)\neq 0$ in $D(s_0,\frac{1}{2})$
for any $y$ is equivalent to  the radius of convergence of the
series (\ref{sec3bis-eq4}) being at least $\frac{1}{2}$ for any
$y$.

Now, the union of the strips $S_1=\{s\in \mathbb{C}:
\frac{1}{2}<\sigma <1 \}$ and $S_2=\{s\in \mathbb{C}:
0<\sigma<\frac{1}{2}  \}$ form the critical strip minus the
critical line $\sigma=\frac{1}{2}$. Moreover, the strip $S=\{s\in
\mathbb{C}: \frac{1}{2}<\sigma <\frac{3}{2} \}$ can be written as

\begin{equation}\label{sec3bis-eq7}
S=\bigcup_{y} D(1+iy,\frac{1}{2})
\end{equation}

We conclude from the above, that if $\zeta(s)$ doest not have a
zero inside the strip $S$ and a fortiori does not have  a zero in
the strip $S_1$, then by the functional equation $\zeta(s)$ cannot
have a zero inside $S_2$ neither. We thus have proved

\begin{theorem}\label{sec3bis-thm1}
The Riemann hypothesis is equivalent to either

\begin{enumerate}
    \item the series (\ref{sec3bis-eq1}) does not have any zero in
    the disk $D(s_0,\frac{1}{2})$ for any $y$.
    \item the radius of convergence of the series (\ref{sec3bis-eq4})
    is
at least $\frac{1}{2}$ for any $y$.
\end{enumerate}
\end{theorem}

\begin{remark}
We have been able to prove that the series (\ref{sec3bis-eq1})
does not have any zero in the disk $D(1,\frac{1}{2})$ (i.e $y=0$).
The proof is trivial and uses the criterion of Petrovitch
\cite{landau} for power series. We also have been able to prove
the well-known result that $D(1,1)$ is a zero-free region. We have
been unable to generalize the proof to any $y$ because as $y$ gets
large the value of $|b_0|$ become very small compared to that of
$|b_1|$. The typical power series and polynomial non-zero regions
criteria are inapplicable. More knowledge on the ratios
$|b_n|/|b_0|$ is needed.
\end{remark}

\begin{remark}
Another  criterion of the Riemann hypothesis can be formulated
using the conformal mapping $s=\frac{1}{1-z}$ which maps the plane
$\Re(s)>\frac{1}{2}$ onto the unit disk $|z|< 1$. We can write

\begin{eqnarray}
(\frac{z}{1-z})\zeta(\frac{1}{1-z})\Gamma(\frac{1}{1-z})\triangleq
f(z) &=&\int_{0}^{\infty}\phi(t)
e^{-t}t^{\frac{z}{1-z}}\,dt\nonumber \\
&=&\sum_{n=0}^{\infty}\tilde{a}_{n}z^{n},\label{sec3bis-eq8}
\end{eqnarray}

where the coefficients $\tilde{a}_n$ are given by

\begin{equation}\label{sec3bis-eq9}
\tilde{a}_n= \int_{0}^{\infty}\phi(t)e^{-t}L_{n}(-\ln{t})\,dt,
\end{equation}

$L_{n}(x)=L_{n}^{0}(x)$ being the Laguerre polynomial of order
$0$.

The Riemann hypothesis is equivalent to the function $f(z)$ having
no zeros in the unit disk.
\end{remark}

Although the above formulations of the Riemann hypothesis seem to
be promising since exact information on the coefficients is known,
we will not pursue this approach. The new approach that we will
adopt is presented next.

\specialsection*{PART II }

In this part we will pursue a completely different approach from
the one presented in PART I. Using the new series representation
of the zeta function of the first part, exact information on its
zeros is provided based on Tauberian-like results.

\section{The Series representation of $(s-1)\zeta(s)$}

In PART I, we showed that $(s-1)\zeta(s)$ process both an integral
and a series representation valid for $\Re(s)>0$. In the remaining
of the paper we will only consider the series representation. We
recall the series representation valid when $\Re(s)>1$:

\begin{equation}\label{sec5-eq1}
\sum_{n=1}^{\infty}\frac{S_{n}(s)}{n+1} =
\frac{s-1}{\Gamma(s)}\int_{0}^{\infty}\frac{e^{-t}t^{s-1}}{1-e^{-t}}\,dt=(s-1)\zeta(s),
\end{equation}

where  $S_{n}(s)$ is given by

\begin{equation}\label{sec5-eq2}
S_{n}(s) = \sum_{k=0}^{n-1}(-1)^{k}{n-1\choose k}(k+1)^{-s}.
\end{equation}

First, we  provide another proof of the validity of the series
representation for $\Re(s)>0$. To prove the analytic continuation
when $\Re(s)>0$, we need to evaluate the sum when $\Re(s)>0$.  The
next lemma, which will also be needed in the rest of the paper,
provides such an estimation. It provides an estimate of the exact
asymptotic order of growth of $\frac{S_{n}(s)}{n+1}$ when $n$ is
large.

\begin{lemma}\label{sec5-lem1}
$\displaystyle{\frac{S_{n}(s)}{n+1}\thicksim \frac{1}{n(n+1)(\log
n)^{1-s}\Gamma(s)}}$ for $n$ large enough and for all
$s=\sigma+it$, $\Re(s)>0$, $ s\notin \{1,2,\cdots\}$.
\end{lemma}
\begin{proof}
By putting $k=m-1$ in (\ref{sec5-eq2}), we have by definition

\begin{eqnarray}
S_{n}(s)&=&\sum_{m=1}^{n}{n-1\choose m-1}(-1)^{m-1}m^{-s}=
\sum_{m=1}^{n}\frac{m}{n}{n\choose m}(-1)^{m-1}m^{-s} \nonumber\\
&=& \frac{-1}{n}\sum_{m=1}^{n}{n\choose m}(-1)^{m}m^{1-s}=
\frac{-1}{n}\Delta_{n}(s-1),\label{sec5-eq3}
\end{eqnarray}

where $\displaystyle{\Delta_{n}(\lambda)\triangleq
\sum_{m=1}^{n}{n\choose m}(-1)^{m}m^{-\lambda}}$.

The asymptotic expansion of sums of the form
$\displaystyle{\Delta_{n}(\lambda)}$, with $\lambda\in\mathbb{C}$
being nonintegral has been given in Theorem 3 of Flajolet et al.
~\cite{flajolet:rice}. With a slight modification of notation, the
authors in ~\cite{flajolet:rice} have shown that
$\Delta_{n}(\lambda)$ has an asymptotic expansion in descending
powers of $\log n$ of the form

\begin{equation}\label{sec5-eq4}
-\Delta_{n}(\lambda)\thicksim (\log
n)^{\lambda}\sum_{j=0}^{\infty}(-1)^{j}
\frac{\Gamma^{(j)}(1)}{j!\Gamma(1+\lambda-j)}\frac{1}{(\log
n)^{j}}
\end{equation}

We apply the theorem to $\Delta_{n}(\lambda)$ with $\lambda=s-1$
to get
\begin{equation}\label{sec5-eq5}
\Delta_{n}(s-1) \thicksim \frac{-(\log{n})^{s-1}}{\Gamma(s)},
\end{equation}

which leads to the result
\begin{equation}\label{sec5-eq6}
S_{n}(s) \thicksim\frac{1}{n(\log n)^{1-s}\Gamma(s)}.
\end{equation}

The Lemma follows from dividing  equation (\ref{sec5-eq6}) by
$n+1$.
\end{proof}

Now to obtain an analytic continuation when $\Re(s)>0$, we simply
observe that the logarithmic test of series in combination with
the asymptotic value of $S_{n}(s)$ provided by
Lemma~\ref{sec5-lem1} imply that the absolute value of the series
on  the left hand side of (\ref{sec5-eq1}) is dominated by a
uniformly convergent series for all finite $s$ whose real part is
greater than 0.

\begin{remark}
By Weierstrass theorem, we can see that the function
$(s-1)\zeta(s)$ can be extended outside of the domain $\Re(s)>1$
and that it does not have any singularity when $\Re(s)>0$.
Moreover, by repeating the same process for $\Re(s)>-k$, $k\in
\mathbb{N}$, it is clear that the series defines an analytic
continuation of $\zeta(s)$ valid for all $s\in\mathbb{C}$.
\end{remark}

\section{Preparation Lemmas}\label{sec6}

Throughout this section, we suppose that $0<\Re(s)<1$. For a fixed
$s=\sigma+it$, we associate with $\zeta(s)$ the following power
series:

\begin{equation}\label{sec6-eq1}
(s-1)\zeta(s,x)\triangleq
\frac{S_{1}(s)}{2}x+\frac{S_{2}(s)}{3}x^{2}+
\cdots+\frac{S_{n-1}(s)}{n}x^{n-1}
+\frac{S_{n}(s)}{n+1}x^{n}+\cdots
\end{equation}
$x\in\mathbb{R}$.

Let's also further define the ``comparison" power series by
\begin{equation}\label{sec6-eq2}
\Phi(x)\triangleq (1-x)\bigg(\frac{\log(1-x)}{-x}\bigg)^{s}=
\phi_{0}+\phi_{1}x+\phi_{2}x^{2}+\cdots+\phi_{n}x^{n}+\cdots
\end{equation}

It is easy to verify that  for  $\sigma>0$

\begin{equation}\label{sec6-eq3}
\lim_{x\to 1} (1-x)\bigg(\frac{\log(1-x)}{-x}\bigg)^{s}=0.
\end{equation}

Furthermore, direct calculation of $\Phi^{\prime}(x)$ yields the
expression

\begin{equation}\label{sec6-eq4}
\Phi^{\prime}(x)=\bigg(\frac{\log(1-x)}{-x}\bigg)^{s}\big(-1+s-\frac{s}{\log(1-x)}-\frac{s}{x}\big).
\end{equation}

Clearly, $\Phi^{\prime}(x)$ is  well-defined for all $x\in [0,1)$
and satisfies

\begin{equation}\label{sec6-eq5}
\lim_{x\to 1}|\Phi^{\prime}(x)| =\infty.
\end{equation}

In other words, the function $\Phi^{\prime}(x)$ is a continuous
well-defined function of $x$, converges for all values of $x\in
[0,1)$ and diverges when $x\to 1$. Moreover, because $\Phi(x)$ is
analytic at $x=0$, $\Phi^{\prime}(x)$ must possess the following
power series expansion around $x=0$:

\begin{equation}\label{sec6-eq6}
\Phi^{\prime}(x)=
\phi_{1}+2\phi_{2}x+\cdots+n\phi_{n}x^{n-1}+\cdots
\end{equation}

Finally, we  associate to the series (\ref{sec6-eq6}) the
following positive coefficients power series:

\begin{equation}\label{sec6-eq7}
\tilde{\Phi}(x)\triangleq
|\phi_{1}|+2|\phi_{2}|x+\cdots+n|\phi_{n}|x^{n-1}+\cdots
\end{equation}

The proofs in the remaining of this section will be based on two
theorems. The first theorem, which is due to  N\"{o}rlund
\cite{norlund} and more recently generalized by Flajolet et al.,
estimates the asymptotic behavior of the coefficients of certain
powers series:

\begin{theorem}[\cite{flajolet:singularities}]\label{sec6-thm1}
Let $\alpha$ be a positive integer and $\beta$ be a real or
complex number, $\beta \notin \{0,1,2,\cdots\}$. Define the
function $f(z)$ by
\begin{equation}\label{sec6-eq8}
f(z)=(1-z)^{\alpha}(\frac{1}{z}\big(\log{\frac{1}{1-z}}\big)^{\beta}.
\end{equation}
Then, the Taylor coefficients $f_{n}$ of $f(z)$ satisfy
\begin{equation}\label{sec6-eq9}
f_{n}\thicksim n^{-\alpha-1}(\log{n})^{\beta}\big(
\frac{e_{1}}{1!}\frac{-(\beta)}{(\log{n})}+
\frac{e_{2}}{2!}\frac{\beta(\beta-1)}{(\log{n})^2}+\cdots\big),
\end{equation}
with
\begin{equation}\label{sec6-eq10}
e_{k}=\frac{d^{k}}{ds^{k}}\bigg(
\frac{1}{\Gamma(-s)}\bigg)\bigg|_{s=\alpha}.
\end{equation}
\end{theorem}

The derivatives in (\ref{sec6-eq10}) when $s=\alpha$ is a positive
integer can be evaluated with the help of the identity:

\begin{equation}\label{sec6-eq10bis}
\frac{1}{\Gamma(-s)}=-\frac{\sin(\pi s)}{\pi}\Gamma(1+s).
\end{equation}

For example, the value of $e_1$ when $\alpha$ is a positive
integer is given by

\begin{equation}\label{sec6-eq10bis2}
e_{1}=-\frac{d}{ds}\bigg( \frac{\sin(\pi
s)}{\pi}\Gamma(1+s)\bigg)\bigg|_{s=\alpha}=-\cos(\pi\alpha)\Gamma(1+\alpha),
\end{equation}

and in this case

\begin{equation}\label{sec6-eq10bis3}
f_{n}\thicksim \frac{\beta
\cos(\pi\alpha)\Gamma(1+\alpha)}{n^{1+\alpha}(\log{n})^{1-\beta}}.
\end{equation}

The second theorem, due to Appell \cite{appell}, is the
counterpart of l'Hospital's rule for divergent positive
coefficients power series:

\begin{theorem}[\cite{borel}~p. 66]\label{sec6-thm2}
Let $f(x), g(x)$ be two  real power series of the form
\begin{equation}\label{sec6-eq11}
f(x)=\sum_{n=1}^{\infty}a_{n}x^{n},
g(x)=\sum_{n=1}^{\infty}b_{n}x^{n}, a_{n}, b_{n} >0 {\rm ~for~
all~} n>N , 0<x<1.
\end{equation}
We further suppose that
\begin{itemize}
 \item  the series
$\sum_{n=1}^{\infty}a_{n}$, $\sum_{n=1}^{\infty}b_{n}$ are both
divergent so that $x=1$ is a singular point of both $f(x)$ and
$g(x)$. \item
$\displaystyle{\lim_{n\to\infty}\frac{a_{n}}{b_{n}}=l}$,
\end{itemize}
then,
\begin{equation}\label{sec6-eq12}
\lim_{x\to 1}\frac{f(x)}{g(x)}=l.
\end{equation}
\end{theorem}

Our first result establishes  an important property on the
behavior of the derivative of the function $\Phi(x)$ when $x$ is
close to $1$:

\begin{lemma}\label{sec6-lem1}
There exists an $x_0\in (0,1)$ and a constant $C$ independent of
$x$ such that for all $x\in (x_0,1)$ we have
$\displaystyle{\frac{|\Phi^{\prime}(x)|}{\tilde{\Phi}(x)}}>C>0$.
\end{lemma}
\begin{proof}
From (\ref{sec6-eq4}),
\begin{equation}\label{sec6-eq13}
|\Phi^{\prime}(x)|=\bigg(\frac{\log(1-x)}{-x}\bigg)^{\sigma}\big|-1+s-\frac{s}{\log(1-x)}-\frac{s}{x}\big|.
\end{equation}

Let's suppose that the  power series expansion of
$\displaystyle{\bigg(\frac{\log(1-x)}{-x}\bigg)^{\sigma}}\triangleq
\Psi(x)$ is given by

\begin{equation}\label{sec6-eq14}
\Psi(x)=
\psi_{0}+\psi_{1}x+\psi_{2}x^{2}+\cdots+\psi_{n}x^{n}+\cdots,
\end{equation}
then applying Theorem~\ref{sec6-thm1} with $\alpha=0$ and
$\beta=\sigma$, implies that for large values of $n$, the
coefficients $\psi_{n}$ satisfy the following asymptotic value:

\begin{equation}\label{sec6-eq15}
\psi_{n}\thicksim \frac{\sigma}{n(\log n)^{1-\sigma}}.
\end{equation}

Similarly, for
$\Phi^{\prime}(x)=\phi_{1}+2\phi_{2}x+\cdots+n\phi_{n}x^{n-1}+\cdots$,
Theorem~\ref{sec6-thm1} with $\alpha=1$ and $\beta=s$, implies
that for  large values of $n$, the coefficients $\phi_{n}$ satisfy
the following asymptotic estimates:

\begin{equation}\label{sec6-eq16}
\phi_{n}\thicksim \frac{-s}{n^{2}(\log n)^{1-s}}.
\end{equation}

The asymptotic value  of $n|\phi_{n}|$ imply by Abel's Theorem and
the logarithmic test of series that, like $\Psi(x)$, the series
$\tilde{\Phi}(x)$ goes to infinity as $x$ approaches 1.

We thus have

\begin{equation}\label{sec6-eq17}
\frac{|\Phi^{\prime}(x)|}{\tilde{\Phi}(x)}
=\frac{\Psi(x)x}{\sum_{n=1}^{\infty}n|\phi_n|x^{n}}
\big|-1+s-\frac{s}{\log(1-x)}-\frac{s}{x}\big|.
\end{equation}

Now as $x$ approaches 1,
$\big|-1+s-\frac{s}{\log(1-x)}-\frac{s}{x}\big|$ approaches 1 so
that given any small $\epsilon>0$ we can find $x_1$ such that for
$x\in (x_1,1)$,
$\big|-1+s-\frac{s}{\log(1-x)}-\frac{s}{x}\big|>1-\epsilon$.

Moreover, since $\Psi(x)x=\sum_{n=1}^{\infty}\psi_{n-1}x^{n}$ and
$\sum_{n=1}^{\infty}n|\phi_n|x^{n}$ both go to infinity as
 $x$ approches 1, and since the asymptotic
estimates (\ref{sec6-eq15})-(\ref{sec6-eq16}) of $\psi_{n-1}$ and
$\phi_n$
 verify
\begin{equation}\label{sec6-eq18}
\lim_{n\to\infty}\frac{\psi_{n-1}}{n|\phi_{n}|}=
\frac{\sigma}{|s|},
\end{equation}
then Theorem~\ref{sec6-thm2} gives

\begin{equation}\label{sec6-eq19}
\lim_{x\to
1}\frac{\sum_{n=1}^{\infty}\psi_{n-1}x^{n}}{\sum_{n=1}^{\infty}n|\phi_{n}|x^{n}}
=\frac{\sigma}{ |s|}.
\end{equation}

In other words, given any small $\epsilon>0$ we can find $x_2$
such that for $x\in (x_2,1)$,
$\frac{\Psi(x)x}{\sum_{n=1}^{\infty}n|\phi_n|x^{n}}>\frac{\sigma}{
 |s|}-\epsilon$.

To complete the proof take for example $C=\frac{\sigma}{ 4 |s|}$
and $x_0=max\{x_1,x_2 \}$.

\end{proof}

The second lemma that we  need establishes a relationship between
the derivative of $\zeta(s,x)$ and that of $\Phi(x)$:

\begin{lemma}\label{sec6-lem2}
Let $\Phi$ be defined as above and let $\zeta^{\prime}(s,x)$ be
$\frac{d}{dx}\zeta(s,x)$, then
\begin{equation}\label{sec6-eq20} \lim_{x\to
1}\frac{(s-1)\zeta^{\prime}(s,x)}{\Phi^{\prime}(x)}=\frac{-1}{
 s\Gamma(s)},
\end{equation}
where the limit is taken from below.
\end{lemma}
\begin{proof}
We have
\begin{equation}\label{sec6-eq21}
(s-1)\zeta^{\prime}(s,x)=\frac{S_{1}(s)}{2}+\frac{2S_{2}(s)}{3}x+
\cdots+\frac{(n-1)S_{n-1}(s)}{n}x^{n-2}
+\frac{nS_{n}(s)}{n+1}x^{n-1}+\cdots
\end{equation}
$x\in\mathbb{R}$.

Lemma~\ref{sec5-lem1} gives
$\displaystyle{\frac{S_{n}(s)}{n+1}\thicksim \frac{1}{n(n+1)(\log
n)^{1-s}\Gamma(s)}}$ for $n$ large enough and for all
$s=\sigma+it$, $0<\Re(s)<1$. Combining with the estimate
(\ref{sec6-eq16}) yields
\begin{equation}\label{sec6-eq22}
\lim_{n\to
\infty}\frac{n\frac{S_{n}(s)}{n+1}}{n\phi_{n}}=\frac{-1}{
s\Gamma(s)}.
\end{equation}

The limit (\ref{sec6-eq22}) is equivalent to saying that there
exists a complex sequence  $\{\epsilon_n\}$ with
$\displaystyle{\lim_{n\to\infty}\epsilon_n=0}$ such that

\begin{equation}\label{sec6-eq23}
n\frac{S_{n}(s)}{n+1}=\frac{-1}{  s\Gamma(s)}n\phi_{n}+ \epsilon_n
n\phi_{n}
\end{equation}
or equivalently

\begin{equation}\label{sec6-eq24}
n\frac{S_{n}(s)}{n+1}x^{n-1} =\frac{-1}{
s\Gamma(s)}n\phi_{n}x^{n-1}+ \epsilon_n n\phi_{n} x^{n-1},
\end{equation}
and finally adding the  equalities for $n=1,2,\cdots$, yields

\begin{equation}\label{sec6-eq25}
(s-1)\zeta^{\prime}(s,x)= \frac{-1}{
s\Gamma(s)}\Phi^{\prime}(x)+\sum_{n=1}^{\infty}\epsilon_n
n\phi_{n} x^{n-1}.
\end{equation}
 Now by dividing all sides of (\ref{sec6-eq25}) by $\displaystyle{
\Phi^{\prime}(x)}$, $0<x<1$, and taking the limit as $x\to 1$, we
get
\begin{equation}\label{sec6-eq26}
\lim_{x\to
1}\frac{(s-1)\zeta^{\prime}(s,x)}{\Phi^{\prime}(x)}=\frac{-1}{
 s\Gamma(s)} +\lim_{x\to
1}\frac{\sum_{n=1}^{\infty}\epsilon_n
n\phi_{n}x^{n-1}}{\Phi^{\prime}(x)}.
\end{equation}

To prove the lemma it suffices  to show that

\begin{equation}\label{sec6-eq27}
\lim_{x\to 1}\bigg|\frac{\sum_{n=1}^{\infty}\epsilon_n
n\phi_{n}x^{n-1}}{\Phi^{\prime}(x)}\bigg|=0.
\end{equation}

Indeed, using our first preparation Lemma~\ref{sec6-lem1}, simple
calculations yield
\begin{eqnarray}
\lim_{x\to 1}\frac{|\sum_{n=1}^{\infty}\epsilon_n
n\phi_{n}x^{n-1}|}{|\Phi^{\prime}(x)|} &\le& \lim_{x\to
1}\frac{\sum_{n=1}^{\infty}n|\epsilon_n|
|\phi_{n}|x^{n-1}}{|\Phi^{\prime}(x)|}\nonumber\\
&\le& \frac{1}{C}\lim_{x\to
1}\frac{\sum_{n=1}^{\infty}n|\epsilon_n
||\phi_{n}|x^{n-1}}{\tilde{\Phi}(x)},\label{sec6-eq28}
\end{eqnarray}

where $\tilde{\Phi}(x)$ is defined in (\ref{sec6-eq7}).

If the series in the numerator is convergent, the result is
obvious. If not, the two series in the right hand side of the last
inequality are both divergent positive coefficients power series.
An application of Theorem~\ref{sec6-thm2} shows that the limit in
(\ref{sec6-eq28}) is equal to the limit of
\begin{equation}\label{sec6-eq29}
\lim_{n\to \infty}\frac{n|\epsilon_n ||\phi_{n}|}{n|\phi_{n}|}=0,
\end{equation}

and the lemma is proved.

\end{proof}

Now let $s$ be a nontrivial zero of $\zeta(s)$. By Abel's theorem
$\lim_{x\to 1}(s-1)\zeta(s,x) = (s-1)\zeta(s,1)=0$. In addition,
$\lim_{x\to 1}\Phi(x)=0$. The next and last preparation lemma
shows that because of the particular function $\Phi(x)$,
l'Hospital's rule which usually does not apply to vector valued or
complex valued function, does apply for this particular case:

\begin{equation}\label{sec6-eq30}
\lim_{x\to 1}\frac{(s-1)\zeta(s,x)}{\Phi(x)}= \lim_{x\to
1}\frac{(s-1)\zeta^{\prime}(s,x)}{\Phi^{\prime}(x)}=\frac{-1}{
 s\Gamma(s)}.
\end{equation}

\begin{lemma}\label{sec6-lem3}
Let $\Phi$ be defined as above and let $s$ be a nontrivial zero of
$\zeta(s)$, then
\begin{equation}\label{sec6-eq31}
\lim_{x\to 1}\frac{(s-1)\zeta(s,x)}{\Phi(x)}=\frac{-1}{
s\Gamma(s)}
\end{equation}
where the limit is taken from below.
\end{lemma}
\begin{proof}
From Lemma~\ref{sec6-lem2}, we have
\begin{equation}\label{sec6-eq32}
\lim_{x\to
1}\frac{(s-1)\zeta^{\prime}(s,x)}{\Phi^{\prime}(x)}=\frac{-1}{
 s\Gamma(s)}.
\end{equation}

Define $\delta(x)$ by

\begin{equation}\label{sec6-eq33}
\delta(x)\triangleq\frac{(s-1)\zeta^{\prime}(s,x)}{\Phi^{\prime}(x)}+\frac{1}{
 s\Gamma(s)}
\end{equation}

so that (\ref{sec6-eq32}) can be written as
\begin{equation}\label{sec6-eq33bis}
\lim_{x\to 1}\delta(x)=0.
\end{equation}

Multiplying equation (\ref{sec6-eq33}) by $\Phi^{\prime}(x)$, and
integrating\footnote{The integral is an improper integral, i.e. it
is defined as $\lim_{\epsilon\to 0}\int_x^{1-\epsilon}f(y)\,dy$}
from $x$ to $1$, we obtain

\begin{equation}\label{sec6-eq34}
(s-1)\big(\lim_{\epsilon\to 0}\zeta(s,1-\epsilon)-\zeta(s,x)\big)+
\frac{\lim_{\epsilon\to 0}\Phi(1-\epsilon)-\Phi(x)}{
s\Gamma(s)}=\int_{x}^{1}\delta(y)\Phi^{\prime}(y)\,dy.
\end{equation}

Now recalling that $\lim_{\epsilon\to 0}\Phi(1-\epsilon)=0$, and
that $s$ is a zero of $\zeta(s)$ so that $\lim_{\epsilon\to
0}\zeta(s,1-\epsilon)=0$, and dividing both sides of
(\ref{sec6-eq34}) by $\Phi(x)$, we finally get
\begin{equation}\label{sec6-eq35}
\frac{(s-1)\zeta(s,x)}{\Phi(x)}+\frac{1}{
s\Gamma(s)}=-\frac{\int_{x}^{1}\delta(y)\Phi^{\prime}(y)\,dy}{\Phi(x)},
\end{equation}

by which we obtain:

\begin{eqnarray}
\lim_{x\to 1}\bigg|\frac{(s-1)\zeta(s,x)}{\Phi(x)}+\frac{1}{
 s\Gamma(s)}\bigg|&\le& \lim_{x\to
1}\frac{\big|\int_{x}^{1}\delta(y)\Phi^{\prime}(y)dy\big|}{|\Phi(x)|}\nonumber\\
&\le&\lim_{x\to
1}\frac{\int_{x}^{1}\big|\delta(y)\big|\big|\Phi^{\prime}(y)\big|\,dy}{|\Phi(x)|}
\label{sec6-eq36}
\end{eqnarray}

By observing that the ratio
 \begin{equation}\label{sec6-eq37}
\frac{\big|\Phi^{\prime}(y)\big|}{\big|\Phi(y)\big|^{\prime}}=
\frac{\big|-1+s-\frac{s}{\log(1-y)}-\frac{s}{y}\big|}{\big|-1+\sigma-\frac{\sigma}{\log(1-y)}-\frac{\sigma}{y}\big|},
\end{equation}
where $\displaystyle{\big|\Phi(y)\big|^{\prime}=
\frac{d}{dy}\big|\Phi(y)\big|}$ is always bounded by a suitable
constant, say  $K$, for $y\in (x_0,1)$, $x_0$ close to 1, the
limit in (\ref{sec6-eq36}) is less than or equal to

\begin{equation}\label{sec6-eq38}
K\lim_{x\to
1}\frac{\int_{x}^{1}\big|\delta(y)\big|\big|\Phi(y)\big|^{\prime}\,dy}{|\Phi(x)|}.
\end{equation}

 The last limit in (\ref{sec6-eq38}) consists of a limit of the ratio
of two real functions that satisfy the hypothesis of l'Hospital's
rule. That is

\begin{equation}\label{sec6-eq39}
 K\lim_{x\to
1}\frac{\int_{x}^{1}\big|\delta(y)\big|\big|\Phi(y)\big|^{\prime}\,dy}{|\Phi(x)|}=K
\lim_{x\to
1}\frac{\big|\delta(x)\big|\big|\Phi(x)\big|^{\prime}}{|\Phi(x)|^{\prime}}=K\lim_{x\to
1}\big|\delta(x)\big|=0.
\end{equation}
 Consequently,

\begin{equation}\label{sec6-eq40}
\lim_{x\to 1}\bigg|\frac{(s-1)\zeta(s,x)}{\Phi(x)}+\frac{1}{
s\Gamma(s)}\bigg|=0,
\end{equation}

and the lemma is proved.
\end{proof}

Equation (\ref{sec6-eq31}) in Lemma~\ref{sec6-lem3} is  quite a
remarquable identity. It says that if $s$ is a nontrivial zero of
$\zeta(s)$ and even though for this particular value of $s$,
$\lim_{x\to 1}(s-1)\zeta(s,x)=0$ and $\lim_{x\to 1}\Phi(x)=0$, the
limit $\displaystyle{\lim_{x\to 1}\frac{
(s-1)\zeta(s,x)}{\Phi(x)}}$ is well-defined and is equal to
$\displaystyle{\frac{-1}{ s\Gamma(s)}}$.

For a nontrivial zero $s$, $\Phi(x)$ is in some sense a measure of
the rate of convergence of $(s-1)\zeta(s,x)$ to zero and hence of
the manner $(s-1)\zeta(s)$ goes to zero. By using the above
identity and comparing the rate of convergence of two symmetric
zeros with respect to the critical line, we can deduce the Riemann
hypothesis.

\section{Proof of Riemann Hypothesis}\label{sec7}

We proceed by contradiction. Suppose that $s=\sigma+it$ is a
nontrivial zero of $\zeta(s)$ with
$\displaystyle{0<\sigma<\frac{1}{2}}$. From the functional
equation (\ref{sec1-eq2}), $1-s$ must also be a nontrivial zero of
$\zeta(s)$.

We have from the previous analysis

\begin{equation}\label{sec7-eq1}
\lim_{x\to 1}\frac{(s-1)\zeta(s,x)}{\Phi(x)}=\frac{-1}{
s\Gamma(s)}, \quad{\rm where}
\end{equation}

\begin{equation}\nonumber
\Phi(x)\triangleq (1-x)\bigg(\frac{\log(1-x)}{-x}\bigg)^{s}.
\end{equation}

Similarly, for $1-s$, we define the comparison function
$\hat{\Phi}(x)$ which is  analogous to $\Phi(x)$ but with $1-s$ in
place of $s$ to get

\begin{equation}\label{sec7-eq2}
\lim_{x\to 1}\frac{-s\zeta(1-s,x)}{\hat{\Phi}(x)}=\frac{-1}{
(1-s)\Gamma(1-s)}, \quad{\rm where}
\end{equation}

\begin{equation}\nonumber
\hat{\Phi}(x)\triangleq (1-x)\big(\frac{\log(1-x)}{-x}\big)^{1-s}.
\end{equation}

Taking absolute values and dividing equation (\ref{sec7-eq1}) by
(\ref{sec7-eq2}) \footnote{We could have evaluated the quotient
directly without taking absolute values, but then the limit of the
quotient exists under the condition that $\lim_{x\to
1}\frac{\hat{\Phi}(x)}{\Phi(x)}$ exists; and this is not always
true for complex-valued functions (e.g. for $s=0.5+it, t\ne 0$,
the limit does not exist). By taking the absolute values, the
quotients are real-valued and we circumvent such cases.}, we must
then have

\begin{equation}\label{sec7-eq3}
\lim_{x\to
1}\bigg|\frac{(1-s)\zeta(s,x)}{s\zeta(1-s,x)}\frac{\hat{\Phi}(x)}{\Phi(x)}\bigg|
=\lim_{x\to
1}\bigg|\frac{(1-s)\zeta(s,x)}{s\zeta(1-s,x)}\bigg|\lim_{x\to
1}\bigg|\frac{\hat{\Phi}(x)}{\Phi(x)}\bigg|=\bigg|\frac{\Gamma(1-s)}{\Gamma(s)}\frac{1-s}{s}\bigg|.
\end{equation}

Now, $\lim_{x\to 1}|\frac{(1-s)\zeta(s,x)}{s\zeta(1-s,x)}|$ is
equal to $|\frac{(1-s)\zeta(s)}{s\zeta(1-s)}|$. The latter
quantity is a finite non-zero value since by continuity with
respect to $s$ the functional equation (\ref{sec1-eq2}) implies
that
\begin{equation}\label{sec7-eq5}
\frac{\zeta(s)}{\zeta(1-s)}=\frac{\pi^{-\frac{1-s}{2}}\Gamma(\frac{1-s}{2})}{\pi^{-\frac{s}{2}}\Gamma(\frac{s}{2})}
\end{equation}

is a finite number. Consequently,

\begin{equation}\label{sec7-eq4}
\lim_{x\to
1}\bigg|\frac{\hat{\Phi}(x)}{\Phi(x)}\bigg|=\bigg|\frac{\zeta(1-s)}{\zeta(s)}\frac{\Gamma(1-s)}{\Gamma(s)}\bigg|
\end{equation}

is a finite number. But exact calculation of the left hand side of
(\ref{sec7-eq4}) gives:

\begin{equation}\label{sec7-eq6}
\lim_{x\to 1}\bigg|\frac{\hat{\Phi}(x)}{\Phi(x)}\bigg|=\lim_{x\to
1}\bigg(\frac{\log(1-x)}{-x}\bigg)^{1-2\sigma}=\infty.
\end{equation}

This contradicts equation (\ref{sec7-eq4}) unless
$\sigma=\frac{1}{2}$ in which case the limit in (\ref{sec7-eq4})
is equal to 1. So there cannot be a zero such that
$0<\sigma<\frac{1}{2}$ and therefore all the zeros must lie on the
line $\sigma=\frac{1}{2}$. The proof is complete.






\end{document}